\documentclass{cimart}

\DeclareMathOperator{\ad}{ad}
\newcommand{\adgo}{\mbox{\rm ad}_{\fg_1}}
\newcommand{\adgn}{\mbox{\rm ad}_{\fg_0}}
\newcommand{\C}{\mathbb{C}}
\newcommand{\cF}{\mathcal{F}}
\newcommand{\cH}{\mathcal{H}}
\newcommand{\cS}{\mathcal{S}}
\newcommand{\fa}{\mathfrak{a}}
\newcommand{\fg}{\mathfrak{g}}
\newcommand{\fh}{\mathfrak{h}}
\newcommand{\fj}{\mathfrak{j}}
\newcommand{\R}{\mathbb{R}}
\newcommand{\sF}{D_2}
\newcommand{\sH}{D_1}
\newcommand{\w}{\omega}

\title{% Please, capitalize only the first word
    K\"ahlerian oxidation and reduction of Lie algebras with  K\"ahler structures
    }

\authors{% Please, use "Firstname Lastname" format, without abreviations
    Mustapha Bachaou, Ignacio Bajo, and Mohamed Louzari
    }

\authorinfo[% Format: please, use "F. Name"
    M. Bachaou]{% Format: please, use "Department (only if 
Abdelmalek Essa\^adi University, %Tetouan, 
    Morocco}{%
    mustapha.bachaou@etu.uae.ac.ma
    }

\authorinfo[% Format: please, use "S. Name"
    I. Bajo]{% Format: please, use "Department (only if 
    University of Vigo, Spain}{%
    ibajo@dma.uvigo.es
    }

\authorinfo[% Format: please, use "T. Name"
    M. Louzari]{% Format: please, use "Department (only if 
    Abdelmalek Essa\^adi University, %Tetouan, 
    Morocco}{%
    m.louzari@uae.ac.ma
    }

\abstract{%
     A K\"ahler structure on a real Lie algebra $\fg$ is a pair $(\omega, J)$ of a 2-cocycle  and a complex structure on $\fg$ which are compatible in the sense that $J$ is skew-symmetric with respect to $\omega$. This paper investigates Kähler structures on Lie algebras, with a particular focus on the nilpotent case. We introduce an algebraic framework for the Kählerian reduction and oxidation of Lie algebras by one-dimensional and two-dimensional ideals. Our main structural result demonstrates that every non-zero quasi-nilpotent Kählerian nilpotent Lie algebra can be inductively reconstructed through a finite sequence of Kählerian central oxidations by lines or planes. This reduction sequence originates from a base algebra that is either $\{0\}$ or one admitting a strongly non-nilpotent complex structure. Moreover, if the complex structure is nilpotent, all intermediate reductions  also have a nilpotent complex structure. 
    }



\acknowledgments{The authors would like to express their sincere gratitude to the referees for the numerous valuable comments and suggestions, which have significantly improved this paper.}

\VOLUME{35}
\YEAR{2027}
\ISSUE{1}
\NUMBER{4}
\DOI{https://doi.org/10.46298/cm.19436}
\editinfo{September 23, 2026}{November 30, 2026}{Ivan Kaygorodov}

\begin{document}

\section*{Introduction}
\addcontentsline{toc}{section}{Introduction}

Kähler Lie algebras provide a natural algebraic framework at the intersection of Lie theory, symplectic geometry, and complex manifolds. A  K\"ahler structure on a real finite-dimensional Lie algebra $\mathfrak{g}$ consists of a symplectic form $\omega$ and an endomorphism $J \in \mathfrak{gl}(\mathfrak{g})$ satisfying $J^2 = -I$, where $I$ stands for the identity map of $\fg$, such that $J$ is integrable, in the sense that $[Jx,Jy]=[x,y]+J[Jx,y]+J[x,Jy]$ holds for all $x,y\in\fg$, and Hamiltonian with respect to the symplectic form,  meaning that $\w(Jx,y)+\w(x,Jy)=0$ for all $x,y\in\fg$. In this paper, we will study K\"ahler structures on Lie algebras and, in particular, on nilpotent Lie algebras. 

A fundamental problem in the study of invariant geometric structures on Lie algebras is to understand how those structures can be constructed via systematic methods. In the particular case of symplectic Lie algebras, Medina and Revoy \cite{MR} and later Dardi\'e and Medina \cite{DM} described the process of symplectic double extension, which provides a method to construct new symplectic $(n+2)$-dimensional Lie algebras from $n$-dimensional ones. These methods of double extension are particular cases of what Baues and Cort\'es called symplectic oxidation in \cite{BC}. The concept of oxidation is the reverse process  of symplectic reduction; one says that a symplectic Lie algebra $\fg_0$ is a reduction of a larger symplectic Lie algebra $\fg$ if there exists an isotropic ideal $\fj$ in $\fg$ such that $\fg_0$ is symplectomorphic to the quotient Lie algebra $\fj^\bot/\fj$. Symplectic oxidation and reduction have already been used to describe  other structures also involving complex structures, such as complex symplectic structures on Lie algebras \cite{Bazz}.

K\"ahlerian structures have been mainly studied in the case where the symmetric bilinear form defined by $g(x,y)=\w (x,Jy)$, usually called the K\"ahler metric, is positive definite. In fact, it is quite common in the literature to reserve the term ``K\"ahler'' for that case and use ``pseudo-K\"ahler'' or ``indefinite K\"ahler'' whenever $g$ can have an arbitrary signature. We will not make such a distinction, so that  we will use the term K\"ahler for both the definite and indefinite cases. Some structure results in the positive definite case were first obtained by Lichnerowicz and Medina \cite{LM} and a complete description of these Lie algebras in terms of double extension was given in \cite{DM2}. However, as far as we know, in the general, possibly indefinite, case, only some particular cases have been studied (see \cite{BBL2, BS, Said, CG, LU} and references therein). 
The main objective of this paper is to introduce and study a systematic framework for K\"ahlerian reduction and oxidation of real Lie algebras, focusing on central oxidations by lines (i.e., one-dimensional vector spaces) and planes (two-dimensional ones) in order to give an inductive description of certain families of K\"ahlerian Lie algebras. 

We pay special attention to the class of nilpotent Lie algebras. For a nilpotent Lie algebra $\fg$ carrying a complex structure $J$, the behavior of $J$ with respect to the central series allows one to distinguish different nilpotency types: strongly non-nilpotent, quasi-nilpotent, and nilpotent complex structures \cite{cord1,LU,LUV1,LUV2}. Explicitly, a complex structure is said to be quasi-nilpotent if the center of the Lie algebra meets its image under the complex structure $J$, and it is said to be strongly non-nilpotent otherwise. Further, associated with the complex structure $J$, one can define the so-called upper $J$-compatible series  $\{\fa_k(J)\}_k$ and one says that the complex structure is nilpotent whenever $\fa_p(J)=\fg$ for some $p\ge 0$ (called the step of nilpotency of $J$) \cite{cord1,LUV1}. Our main structural result establishes that every non-zero quasi-nilpotent K\"ahlerian nilpotent Lie algebra  can be reconstructed through a finite sequence of K\"ahlerian central oxidations by lines or planes, starting from a base algebra which is either $\{0\}$  or a K\"ahlerian Lie algebra with a strongly non-nilpotent complex structure. Furthermore, when $J$ is nilpotent, all intermediate reductions preserve the nilpotency of $J$ and the process necessarily terminates at $\{0\}$.

The paper is organized into four sections following this introduction. 
In Section 1, we gather essential preliminaries on symplectic Lie algebras, left-symmetric products, complex structures, and K\"ahler structures. In particular, we recall the notions of symplectic reduction and oxidation and that of symplectic double extension.
    Section 2 is devoted to examining central symplectic oxidations by two-dimensional planes and to characterizing the set of admissible data $(D_1, D_2, (a_{ij}))$ that yields a well-defined extended Lie algebra.
 In Section 3, we formalize K\"ahlerian reductions and oxidations by lines and planes, detailing the exact compatibility equations required for the extended complex structure to remain integrable and skew-symmetric with respect to the symplectic form.
   Finally, in Section 4, we focus on K\"ahlerian nilpotent Lie algebras. We construct explicit examples of strongly and weakly non-nilpotent K\"ahler structures, prove the main inductive decomposition theorem (Theorem \ref{tr:mainth}), and analyze the case of $2$-step nilpotent Lie algebras, characterizing the conditions under which the step of the complex structure is $\nu(J) \in \{2, 3\}$.

\section{Preliminaries}

Throughout this paper, we consider real Lie algebras of finite dimension. For a Lie algebra $\fg$, we will denote by $Z(\fg)$ its center and $\mbox{Der}(\fg)$ the Lie algebra of its derivations.

\begin{definition}{ A {\it symplectic form} on a Lie algebra $\fg$ is a non-degenerate skew-symmetric bilinear form $\omega:\fg\times\fg\to\R$ which is {\it closed} in the sense that
$$\omega([x,y],z)+\omega([y,z],x)+\omega([z,x],y)=0$$
holds for all $x,y,z\in\fg$. Equivalently, $\omega$ is a non-degenerate 2-cocycle for the Chevalley-Eilenberg scalar cohomology of $\fg$. The pair $(\fg,\w)$ is called a {\it symplectic Lie algebra}.

In this case, if $F\in {\mathfrak{gl}}(\fg)$, we will denote by $F^*\in {\mathfrak{gl}}(\fg)$ its adjoint map with respect to the form $\w$, uniquely defined by the identities $\w(F^*x,y)=\omega(x,Fy)$ for every $x,y\in\fg$. Therefore, $F$ will be called {\it Hamiltonian} if $F+F^*=0$ and {\it skew-Hamiltonian} if $F=F^*$. It is clear that, if $F\in {\mathfrak{gl}}(\fg)$ is an arbitrary map, the linear endomorphism $\cS(F)=F+F^*$ is skew-Hamiltonian.

Let $(\fg_1,\omega_1)$, $(\fg_2,\omega_2)$ be two symplectic Lie algebras. We say that $\Phi:\fg_1\to\fg_2$ is a {\it symplectomorphism} if  $\Phi$ is an isomorphism of Lie algebras such that $$\omega_2(\Phi x,\Phi y)=\omega_1(x,y)$$ for all $x,y\in\fg_1$.}

\end{definition}

\begin{remark} { On every symplectic Lie algebra $(\fg,\w)$, one can define a new product \newline $\cdot:\fg\times\fg\to\fg$ by the identity
$$\w(x\cdot y,z)=-\w(y,[x,z]),\quad x,y,z\in\fg.$$
It is well known that such a product is left-symmetric, in the sense that 
$$[x,y]\cdot z=x\cdot(y\cdot z)-y\cdot (x\cdot z)$$
is verified for every $x,y,z\in\fg$. Further, a simple computation shows that $[x,y]=x\cdot y-y\cdot x$, so that $\fg$ admits a {\it pre-Lie algebra} structure in the sense of \cite{burde} and references therein. Notice that, geometrically speaking, this is equivalent to saying that each associated Lie group can be endowed with the left-invariant flat and torsion-free connection defined by $\nabla^\w_xy=x\cdot y$.
}
\end{remark}

In \cite{BC}, the concepts of symplectic reduction and symplectic oxidation were introduced. 

\begin{definition}{ 
Let $(\fg,\omega)$ be a symplectic Lie algebra and let $\fj \subseteq \fg$ be an isotropic ideal. Its symplectic orthogonal
\[
\fj^{\perp}=\{x\in\fg \mid \omega(x,\fj)=\{0\}\}
\]
is a Lie subalgebra of $\fg$ containing $\fj$. Consequently, the form $\omega$ induces a symplectic form $\bar{\omega}$ on the quotient Lie algebra
$
\bar{\fg}=\fj^{\perp}/\fj .
$
The symplectic Lie algebra $(\bar{\fg},\bar{\omega})$ is called the {\it symplectic reduction} of $(\fg,\omega)$ with respect to the isotropic ideal $\fj$.

If $\fj$ is a central ideal, the reduction is then called \emph{central}. In such a  case, $\fj^{\perp}$ contains the commutator ideal $[\fg,\fg]$, so that it is also an ideal of $\fg$.}
\end{definition}

\begin{remark} { Not every symplectic Lie algebra admits a reduction (see, for instance, \cite[Section 2.4]{BC}). However, if $(\fg,\omega)$ admits a one-dimensional ideal, then it is automatically isotropic and hence $\fg$ can be reduced. This is the case of symplectic Lie algebras with non-trivial center. Therefore, every nilpotent Lie algebra is reducible.  
}
\end{remark}

\begin{definition} { One says that a symplectic Lie algebra $(\fg,\omega)$ is a {\it symplectic oxidation} of another symplectic Lie algebra $(\fg_0,\w_0)$ if and only if $\fg$ admits a reduction $(\bar\fg,\bar\w)$ symplectomorphic to $(\fg_0, \w_0)$. If the reduction $\bar\fg$ is central, we say that  $(\fg,\omega)$ is a {\it central symplectic oxidation} of $(\fg_0,\w_0)$.}
\end{definition}

\begin{remark} \label{nota3} { Notice that non-symplectomorphic symplectic Lie algebras can have symplectomorphic reductions. For instance, the Lie algebra ${\mathcal K}=\mbox{$\R$-span}\{x_1,x_2,x_3,x_4\}$ with non-trivial bracket $[x_1,x_2]=x_3$ endowed with the skew-symmetric bilinear form $$\omega_{\mathcal K}( x_1,x_4)=\omega_{\mathcal K}( x_2,x_3)= 1$$ is symplectic and its reduction with respect to the central ideal $\fj=\R x_4$ is the abelian Lie algebra $\R^2$ with its unique-up-to-scalar symplectic form. But this reduction is the same as the one obtained from the abelian Lie algebra $\R^4$ furnished with any symplectic form when one considers $(\R a)^\bot/\R a$ for any vector $a\in\R^4$.

This obviously means that a symplectic Lie algebra $(\fg_0,\w_0)$ admits, as a rule, several different oxidations. 
}
\end{remark}

When  $(\fg_0, \w_0)$ is a symplectic reduction of $(\fg,\w)$ with respect to a one-dimensional central ideal, we say that $(\fg,\w)$ is a {\it one-dimensional central oxidation} of $(\fg_0,\w_0)$ \cite[Section~1.4]{BC}. This case was already described by Dardi\'e and Medina in \cite{DM} under the name of {\it symplectic double extension}, and can be characterized by the following result which is just a reformulation of the results in \cite{BC,DM}:

\begin{proposition}\label{sde} Let $(\fg_0,[\ ,\ ]_0,\w_0)$ be a symplectic Lie algebra and $(D,a)\in \mbox{\rm Der}(\fg_0)\times\fg_0$ be such that $\cS(D)D-\ad_{\fg_0}(a)$ is Hamiltonian. 
On the vector space  $\fg={\mathbb R}z\oplus\fg_0\oplus {\mathbb R}v$ define a skew-symmetric bracket  by
$$[z,v]=[z,x]=0\, ,\quad[x,y]=[x,y]_0-\omega_0 (\cS(D)x,y)z\, ,\quad [v,x]=Dx+\omega_0 (x,a)z,\quad x,y\in\fg_0,$$
 and consider the skew-symmetric bilinear form $\omega$ defined by
$$ \omega (z,v)=1\, ,\quad \omega (x,z)= \omega (x,v)=0\, ,\quad\omega (x,y)=\omega_0 (x,y) ,\quad  x,y\in\fg_0.$$
Then $(\fg,\omega)$ is a symplectic Lie algebra whose reduction with respect to $\fj=\R z$ is $(\fg_0,\w_0)$.
Moreover, every one-dimensional central oxidation of $(\fg_0,\w_0)$ is constructed in that way.
\end{proposition}

\begin{remark} { Notice that the condition $\cS(D)D-\ad_{\fg_0}(a)$ Hamiltonian is equivalent to 
$$\w_0(((D+D^*)D+D^*(D+D^*))x,y)=\w_0([a,x]_0,y)+\w_0(x,[a,y]_0)=\w_0(a,[x,y]_0),\,\, x,y\in\fg_0,$$
and this means that the 2-form $\Omega(x,y)=\w_0(((D+D^*)D+D^*(D+D^*))x,y)$ is a 2-coboundary.

In \cite[Def. 5.1]{BB-PK} the construction above was generalized to give a description of one-dimensional oxidations which are not necessarily central.}
\end{remark}

\begin{definition} { Let $\fg$ be a real Lie algebra and denote by $I$ its identity map. A linear map $J\in {\mathfrak{gl}}(\fg)$ 
 such that $J^2=-I$ and
$[Jx,Jy]=[x,y]+J[Jx,y]+J[x,Jy]$ for all $x,y\in\fg$
is said to be a { complex structure} on $\fg$. 
}
\end{definition}

\begin{remark}{ Two particular cases of complex structures have been thoroughly studied. These are {\it bi-invariant} complex structures, which are those verifying $J[x,y]=[Jx,y]$ for all $x,y\in\fg$, and {\it abelian} complex structures, which verify $[Jx,Jy]=[x,y]$ for all $x,y\in\fg$ \cite{Andrada1,ABD,BS,DF}. It is well known that there exists a bi-invariant complex structure on a Lie algebra if and only if the algebra admits a structure of complex Lie algebra. Otherwise, one easily sees that a complex structure $J$ on a  Lie algebra $\fg$ is abelian if and only if the eigenspaces of  $J$  are abelian subalgebras of the complexification $\fg^\C$; this forces the Lie algebra $\fg$ to be 2-step solvable.}\end{remark}

We will mainly work with complex structures on nilpotent Lie algebras. It is known that a complex structure need not preserve the ideals of  the  central descending (or ascending) series. Thus, in order to study different behaviors,  some different classes of complex structures were distinguished  in \cite{LUV1}:

\begin{definition}{
Let $J$ be a complex structure on a nilpotent Lie algebra $\fg$.
The \textit{upper $J$-compatible series} of $\fg$ is the sequence $ \{ \fa_k(J) \}_k $ defined by
\begin{equation}\label{ascending_J-compatible_series}
\fa_0(J) = \{0\}, \quad
\fa_k(J) = \{ X \in \fg \mid [X, \fg] \subset \fa_{k-1}(J) ,\ [JX, \fg] \subset \fa_{k-1}(J) \}, 
\quad k \geq 1.
\end{equation}
It is not difficult to see that each $\fa_k(J)$ is an ideal and that $\fa_{k}(J)\subset\fa_{k+1}(J)$ for all $k\ge 0$, although the inclusion need not be strict \cite{cord1}.

The complex structure $J$ is called \textit{strongly non-nilpotent} if \( \fa_1(J) = \{0\} \) and 
    \textit{quasi-nilpotent} otherwise. Furthermore, $J$ is said to be
     \textit{nilpotent} if there exists an integer \( p > 0 \) such that $\fa_p(J) = \fg$. In this last case, the smallest integer verifying such identity is sometimes called the {\it step of the complex structure},  and following \cite{Gao} we will denote it by $\nu(J)$. A quasi-nilpotent complex structure that is not nilpotent is called {\it weakly non-nilpotent} \cite{LUV1}.
}    
\end{definition}

\begin{remark}{ Note that $\fa_1(J)=Z(\fg)\cap JZ(\fg)$ and, as a consequence, bi-invariant or abelian complex structures on nilpotent Lie algebras are always quasi-nilpotent. In fact, both types of complex structures on nilpotent Lie algebras are always nilpotent because  $[JX,\fg]=\pm[X,J\fg]=[X,\fg]$ for all $X\in\fg$ implies that the upper $J$-compatible series coincides with the central ascending series of the Lie algebra.}
\end{remark}

\begin{definition}{ 
Let $\fg$ be a  Lie algebra, $J$ a complex structure  and $\w$ a symplectic form  on $\fg$. If 
$\w(Jx,y)=-\w(x,Jy)$  for all $x,y \in \fg$, then $(\fg, J, \w)$ is a { K\"ahler Lie algebra}. 

Equivalently, if $g(x,y)=\w(x,Jy)$, the triple $(\fg,J,g)$ is K\"ahler if $J$ is $g$-orthogonal and parallel with respect to the Levi-Civita connection defined by $g$. This is to say, if for any $x,y\in\fg$ we define $\nabla_xy\in\fg$ by the Koszul formula
$$2g(\nabla_xy,z)=g([x,y],z)-g([y,z],x)+g([z,x],y),\quad z\in\fg,$$ then
$g(Jx,Jy)=g(x,y)$ and $\nabla_xJy=J\nabla_xy$ for all $x,y\in\fg$. The bilinear form $g$ is usually called the {\it K\"ahler metric}.

In the literature, the term ``K\"ahler'' is frequently restricted to the case where the metric $g$  is positive definite and ``pseudo-K\"ahler'' is used for the (possibly) indefinite cases. In this work, we will not make a distinction, so that our K\"ahler metrics can be indefinite.

If the  nilpotent Lie algebra $\fg$ has a K\"ahlerian structure $(J,\w)$, we will say that the K\"ahlerian structure is {\it nilpotent/quasi-nilpotent/strongly non-nilpotent} according to the nilpotency type of $J$.}
\end{definition}

It is quite obvious that if $(\fg_1,J_1,\w_1)$ and $(\fg_2,J_2,\w_2)$ are two K\"ahlerian Lie algebras, then the direct product algebra $\fg_1\oplus\fg_2$ will also be K\"ahlerian when furnished with  the  pair $(J,\w)$ defined by $J_{|\fg_i}=J_i$, $\w_{|\fg_i\times\fg_i}=\w_i$ and $\w(\fg_1,\fg_2)=\{0\}$. This justifies the following definition:

\begin{definition} { A K\"ahlerian Lie algebra $(\fg,J,\w)$ is said to be  {\it K\"ahlerian decomposable} if it decomposes as the direct sum  $\fg=\fh_1\oplus\fh_2$ of two ideals $\fh_1$, $\fh_2$ which are non-degenerate for $\w$ and stable under $J$. Note that this is equivalent to the existence of a  proper non-zero ideal $\fh$ of $\fg$ which is $J$-stable, non-degenerate with respect to $\w$ and whose $\w$-orthogonal is also an ideal.

A K\"ahlerian Lie algebra that is not decomposable is said to be (K\"ahlerian) {\it indecomposable}.}
\end{definition} 

\begin{remark}{ If $(\fg,J,\w)$ is K\"ahler and $J$ is bi-invariant, then $\fg$ must be abelian. This is not the case when $J$ is abelian and there are large families of K\"ahler Lie algebras with abelian complex structure \cite{BS,BS-AGAG,BS-TG}.}\end{remark}

The aim of this paper is the study of reduction and oxidation methods for the construction and description of (nilpotent) Lie algebras with K\"ahler structures. It is convenient to start with the study of symplectic oxidations by 2-dimensional central ideals, which we will call {\it central symplectic oxidations by planes}. 

\section{Central symplectic oxidations by planes}

%\noindent {\it Notation.} In the sequel, if $(\fg_0,\w_0)$ is a symplectic Lie algebra, for $(D,a)\in {\mathfrak{gl}}(\fg_0)\times\fg_0$, we will denote
%$$\cH(D,a)=\cS(D)-\ad_{\fg_0}(a)=(D+D^*)D-\ad_{\fg_0}(a).$$
 %
%

 \begin{definition}\label{compd}{ 
Let $(\fg_0,\omega_0)$ be a symplectic Lie algebra. We say that $$(\sH,\sF, a_{11},a_{12},a_{21},a_{22})\in \mbox{\rm Der}(\fg_0)^2\times\fg_0^4$$ is a set of {\it admissible data} for a central symplectic oxidation by a plane if the following conditions hold:
\begin{enumerate}
\item[(i)]  $\cS(D_i)D_j-\ad_{\fg_0}(a_{ij})$ is Hamiltonian for $i\le j$,
\item[(ii)] $[D_1,D_2]=\ad_{\fg_0}(a_{12}-a_{21})$,
\item[(iii)] $D_1^*a_{22}=D_2^*a_{21}-\cS(\sF)(a_{12}-a_{21}) \text{ and } D_2^*a_{11}=D_1^*a_{12}+\cS(\sH)(a_{12}-a_{21})$.
\end{enumerate}
}
\end{definition}

\begin{remark}\label{D21}{ Note that conditions (i) and (ii) above  imply that $\cS(D_2)D_1- \ad_{\fg_0}(a_{21})$ is also Hamiltonian. Actually, if $\cH=\cS(D_2)D_1- \ad_{\fg_0}(a_{21})$,  one has
\begin{align*}
    \cH+\cH^*&=\sF\sH+\sF^*\sH-\ad_{\fg_0}(a_{21})+\sH^*\sF^*+\sH^*\sF-\ad_{\fg_0}(a_{21})^*\\
&=  \sF\sH+\sF^*\sH+[\sH,\sF]-\ad_{\fg_0}(a_{12})+\sH^*\sF^*+\sH^*\sF+[\sF^*,\sH^*]-\ad_{\fg_0}(a_{12})^*\\
&= \sH\sF-\ad_{\fg_0}(a_{12})+\sF^*\sH+\sH^*\sF+\sF^*\sH^*-\ad_{\fg_0}(a_{12})^*\\
&= (\cS(\sH)\sF-\ad_{\fg_0}(a_{12}))+(\cS(\sH)\sF-\ad_{\fg_0}(a_{12}))^*=0.
\end{align*}}
\end{remark}

\begin{proposition} \label{P1}
Let \( (\fg_0,[\ ,\ ]_0, \omega_0) \) be a symplectic Lie algebra, and let $$(\sH,\sF, a_{11},a_{12},a_{21},a_{22})\in \mbox{\rm Der}(\fg_0)^2\times\fg_0^4$$ be a set of admissible data for a central symplectic oxidation by a plane.
For $\alpha,\beta\in\R$, let us define on the vector space $
\fg = \mathbb{R}z \oplus \mathbb{R}z' \oplus \fg_0 \oplus \mathbb{R}v' \oplus \mathbb{R}v$
a bilinear map \( [\cdot , \cdot]: \fg \times \fg \to \fg \) by the following possibly non-zero brackets:
\begin{align*}
[v, v'] &= \alpha z + \beta z' + a_{12}-a_{21}, \\
[v, x] &=   \sH x+\omega_0(x,a_{11})z + \omega_0(x,a_{21})z', \\
[v', x] &=\sF x+ \omega_0(x,a_{12})z + \omega_0(x,a_{22})z', \\
[x,y] &= [x,y]_0 - \omega_0(\cS(\sH)x,y)z - \omega_0(\cS(\sF)x,y)z',
\end{align*}
 for all \( x,y \in \fg_0 \), and let \( \omega: \fg \times \fg \to \mathbb{R} \) be the skew-symmetric form defined by the non-vanishing couplings:
\begin{eqnarray*}
\omega(z, v) &= \omega(z', v') = 1, \quad 
\omega(x, y) &= \omega_0(x, y), \quad  x, y \in \fg_0.
\end{eqnarray*}
Then \( (\fg, [\cdot, \cdot],\omega) \) is a symplectic Lie algebra. 
Moreover,  $(\fg_0,[\cdot, \cdot]_0,\omega_0)$ is the reduction of \( (\fg, [\cdot, \cdot], \omega) \) with respect to the central ideal $\fj=\R z\oplus\R z'$.
\end{proposition}

\begin{proof}
 Denote by $(\fg_1,[\ ,\ ]_1,\omega_1)$ the central oxidation of $(\fg_0,\omega_0)$ by means of $(\sF,a_{22})$ constructed as in Proposition \ref{sde}. One then has $\fg_1=\R z'\oplus\fg_0\oplus\R v'$ with  the nonzero brackets and pairings of $\omega_1$ given by
\begin{align}
& [x,y]_1=[x,y]_0- \omega_0(\cS(\sF)x, y)z', \quad  [v',x]_1=\sF x+ \omega_0(x,a_{22}) z',\label{corchg1}\\ 
& \omega_1(x,y)=\omega_0(x,y),\quad \omega_1(z',v')=1,\nonumber
\end{align}
for $x,y\in\fg_0$.
 Define a linear map $\cF_1:\fg_1 \to \fg_1$  by
$$
\cF_1(z')=0, \quad
\cF_1(v') = \beta z' + a_{12} - a_{21}, \quad
\cF_1(x) = \sH x + \omega_0(x,a_{21})z',$$
for $ x \in \fg_0$.

Let us see that $\cF_1$ is a derivation of $\fg_1$. One clearly has $\cF_1[z',X]=[\cF_1z',X]+[z',\cF_1X]$ for all $X\in\fg_1$ and for every $x\in\fg_0$ we have that
\begin{align*}
    & \cF_1[v',x]_1-[\cF_1v',x]_1-[v',\cF_1x]_1 \\ 
& =[\sH,\sF](x) - \adgn(a_{12}-a_{21})(x) +\omega_0(x, \sF^*a_{21} -(\cS(\sF)(a_{12}-a_{21}) + \sH^*a_{22}))z'= 0,
\end{align*}
because $\sF^*a_{21} = \cS(\sF)(a_{12}-a_{21}) + \sH^*a_{22}$ and $[\sH,\sF]=\adgn(a_{12}-a_{21})$, as in Definition~\ref{compd}. Besides, using that $\sH$ is a derivation of $\fg_0$ and Remark \ref{D21}, we  have for all $x,y\in\fg_0$ that
\begin{eqnarray*}
&& \cF_1 [x,y]_1 - [\cF_1 x,y]_1- [x,\cF_1 y]_1\\
&& =\sH [x,y]_0 - [\sH x,y]_0- [x,\sH y]_0  +  \omega_0([x,y]_0,a_{21})z' +\omega_0(\cS(\cS(\sF) \sH) x,y)z'\\
&& =-  \omega_0([y,a_{21}]_0,x)z' -  \omega_0([a_{21},x]_0,y) z'+\omega_0(\cS(\cS(\sF) \sH) x,y)z' \\&&=\omega_0(\cS(\cS(\sF) \sH-\adgn(a_{21})) x,y)z'= 0.
\end{eqnarray*}
Thus, $\cF_1$ is a derivation of $\fg_1$, as claimed. 

A straightforward calculation shows that the $\omega_1$-adjoint of $\cF_1$ is defined by
$$\cF_1^{*}(z')=0 , \quad
\cF_1^{*}(v') = -\beta z' + a_{21}, \quad
\cF_1^{*}(x) = \sH^{*} x + \omega_0(x,a_{12}-a_{21})z',$$
for every $x\in\fg_0$.
Therefore $(\cS(\cF_1)\cF_1+\cF_1^*\cS(\cF_1)) z'=0$ and, using the conditions given in Definition \ref{compd}, a direct calculation gives
\begin{eqnarray*}
&& (\cS(\cF_1)\cF_1+\cF_1^*\cS(\cF_1)) v'=\cS(\sH)(a_{12}-a_{21})+\sH^*a_{12}=\sF^*a_{11},\\
&&(\cS(\cF_1)\cF_1+\cF_1^*\cS(\cF_1)) x=\cS(\cS(\sH)\sH)x+\omega_0(x,\cS(\sH)(a_{12}-a_{21})+\sH^*a_{12})z'\\& & =\cS(\adgn(a_{11}))x+\omega_0(x,\sF^*a_{11})z'.
 \end{eqnarray*}
But if one considers $b_1=a_{11}-\alpha z'$, one gets
\begin{eqnarray*}
& \adgo(b_1)z'=0,\,\,\, \adgo(b_1)v'=-\sF a_{11}-\omega_0(a_{11},a_{22})z',\,\,\,&\\
& \adgo(b_1)x=\adgn(a_{11})x-\omega_0(\cS(\sF)a_{11},x)z',\,\,\,\adgo(b_1)^*z'=0,\,\,\,&\\
&  \adgo(b_1)^*v'=\cS(\sF)a_{11}+\omega_0(a_{11},a_{22})z',\,\,\, \adgo(b_1)^*x=\adgn(a_{11})^*x+\omega_0(\sF a_{11},x)z',&\end{eqnarray*}
from which one deduces
$$\cS(\cF_1)\cF_1+\cF_1^*\cS(\cF_1)=\adgo(b_1)+\adgo(b_1)^*.$$
This means that $\cS(\cF_1)\cF_1-\adgo(b_1)$ is Hamiltonian, and we can construct the oxidation of $(\fg_1,\w_1)$ by means of $(\cF_1,b_1)$ as in Proposition \ref{sde}. The extended Lie algebra will then be
$\fg=\R z \oplus  \fg_1 \oplus  \mathbb{R}v$ with the bracket
$$[z,v]=[z,X]=0\, ,\quad[X,Y]=[X,Y]_1-\omega_1 (\cS(\cF_1)X,Y)z\, ,\quad [v,X]=\cF_1X+\omega_1 (X,b_1)z,$$ for $X,Y\in\fg_1.$
It suffices to use that $\fg_1=\R z' \oplus  \fg_0 \oplus  \mathbb{R}v'$ and (\ref{corchg1}) to see that $(\fg,[\ ,\ ],\omega)$ is exactly as claimed in the statement, which is, as a consequence, a symplectic Lie algebra.
The fact that $(\fg_0,\w_0)$ is the reduction of $(\fg,\w)$ with respect to the central ideal $\fj=\R z\oplus\R z'$ is now obvious since $\fj^{\perp} =\R z\oplus\R z'\oplus\fg_0$.
\end{proof}

\begin{example}\label{ejem1} { Notice that distinct choices of compatible data (even with the same derivations $D_1,D_2$ and the same values of $\alpha,\beta\in\R$) can give rise to non-isomorphic Lie algebras, as we  show next. Let us consider $\fg_0=\R^2=\R x_1\oplus\R x_2$ with the symplectic form $\omega_0(x_1,x_2)=1$ and take $D_1=D_2=0$, so that we get a set of compatible data for every $(a_{11},a_{12},a_{21},a_{22})\in\fg_0^4.$ Let us consider $\alpha=\beta=0$. Choosing $a_{ij}=0$ for all $i,j\le 2$, one clearly obtains the abelian Lie algebra $\R^6$. If we choose $a_{11}=-x_1$, $a_{12}=x_2$, $a_{21}=a_{22}=0$, then the non-trivial brackets are $[v,v']=x_2$, $[v,x_2]=[v',x_1]=z$, so that  the Lie algebra is a 3-step nilpotent symplectic Lie algebra. For $a_{11}=-x_1$, $a_{22}=x_2$, $a_{21}=a_{12}=0$, we have the non-trivial brackets $[v,x_2]=z$, $[v',x_1]=z'$, so that the Lie algebra is isomorphic to the direct sum of two 3-dimensional Heisenberg algebras ${\mathfrak h}_3\oplus{\mathfrak h}_3$.
}\end{example}

\begin{proposition}\label{P2}
Every central symplectic oxidation by a plane of  a symplectic Lie algebra \( (\fg_0, \omega_0) \) is constructed as in Proposition \ref{P1} for some $\alpha,\beta\in\R$ and certain admissible data $(D_1,D_2, (a_{ij})_{i,j\le 2})\in \mbox{\rm Der}(\fg_0)^2\times\fg_0^4$.
\end{proposition}

\begin{proof}
Let $(\fg,[\ ,\ ],\omega)$ be a central symplectic oxidation by a plane of  \( (\fg_0, [\ ,\ ]_0,\omega_0) \). Then, $\fg_0$ is the symplectic reduction of $\fg$ with respect to some  ideal $\fj=\R z\oplus\R z'$, where $z, z' \in Z(\fg)$ verify $\omega(z,z')=0.$ Since $\w$ is non-degenerate, we can find  $v,v'\in\fg$ such that $\w(z,v)=\w(z',v')=1$ and $\w(z,v')=\w(z',v)=\w(v,v')=0$. The vector space $W=\R z\oplus\R z'\oplus \R v'\oplus\R v$ is non-degenerate with respect to $\w$, so that  its symplectic orthogonal $W^\bot$ must also be non-degenerate. Notice that $\fj^\bot=\fj\oplus W^\bot$, so that $\fg_0=\fj^\bot/\fj$ is symplectically isomorphic (as a vector space) to $W^\bot$ and, hence, we can transfer the Lie structure of $\fg_0$ to $W^\bot$. This means that  we have the vector space decomposition $\fg=\R z\oplus\R z'\oplus\fg_0\oplus\R v'\oplus\R v$ and the symplectic form is given by the non-zero couplings:
$$\w(z,v)=\w(z',v')=1,\quad \w(x,y)=\w_0(x,y),$$
for $x,y\in\fg_0$. It is clear that $[\fg,\fg]\subset Z(\fg)^\bot\subset \fj^\bot=\fj\oplus \fg_0$. But this means that there exist $x_0\in\fg_0$, $\alpha,\beta\in\R$, $D_1,D_2\in {\mathfrak{gl}}(\fg_0)$, linear forms $f,g,f',g':\fg_0\to\R$,  and skew-symmetric bilinear forms $\varphi_1,\varphi_2:\fg_0\times\fg_0\to\R$ such that the {\it a priori} non-trivial brackets in $\fg$ will be
\begin{eqnarray*}
& & [v,v']=x_0+\alpha z+\beta z',\quad [v,x]=D_1x+f(x)z+g(x)z', \\
& & [v',x]=D_2x+f'(x)z+g'(x)z', \quad [x,y]=[x,y]_0+\varphi_1(x,y)z+\varphi_2(x,y)z'.
\end{eqnarray*}
A straightforward calculation shows that $[[v,v'],x]+[[v',x],v]+[[x,v],v']=0$ with $x\in\fg_0$ gives $[D_1,D_2]=\adgn(x_0)$ and
\begin{equation}
\label{jacobi1}
 \varphi_1(x_0,x)=-f'(D_1x)+f(D_2x),\quad \varphi_2(x_0,x)=-g'(D_1x)+g(D_2x).
\end{equation}
Similarly, the identities $$[[v,x],y]+[[x,y],v]+[[y,v],x]=[[v',x],y]+[[x,y],v']+[[y,v'],x]=0$$ give that $D_1,D_2$ are derivations of $\fg_0$ and that 
\begin{eqnarray}
\label{jacobi2}
&& \varphi_1(D_1x,y)+\varphi_1(x,D_1y)=f([x,y]_0), \quad \varphi_2(D_1x,y)+\varphi_2(x,D_1y)=g([x,y]_0)\\
\label{jacobi3}
&& \varphi_1(D_2x,y)+\varphi_1(x,D_2y)=f'([x,y]_0), \quad \varphi_2(D_2x,y)+\varphi_2(x,D_2y)=g'([x,y]_0).
\end{eqnarray}

Besides, since $\w$ is a symplectic form, the condition
$$\w([v_1,v_2],v_3)+\w([v_2,v_3],v_1)+\w([v_3,v_1],v_2)=0$$ should be verified for all $v_1,v_2,v_3\in\fg$.
When we take $v_i\in\fg_0$ for all $i\le 3$, we just get that $\w_0$ is also symplectic, while the cases $(v_1,v_2,v_3)=(v,x,y)$,  $(v_1,v_2,v_3)=(v',x,y)$, $(v_1,v_2,v_3)=(v,v',x)$, with $x,y\in\fg_0$, are respectively equivalent to
\begin{eqnarray} 
& & \varphi_1(x,y)=-\w_0( \cS(D_1)x,y),\quad \varphi_2(x,y)=-\w_0( \cS(D_2)x,y),\label{sympl1}\\ & &  \w_0(x_0,x)=f'(x)-g(x).\label{sympl2}
\end{eqnarray}

As $\w_0$ is non-degenerate, there exist $a_{11},a_{12},a_{21},a_{22}$ such that
$$f(x)=\w_0(x,a_{11}),\, f'(x)=\w_0(x,a_{12}),\, g(x)=\w_0(x,a_{21}),\, g'(x)=\w_0(x,a_{22}).$$
With such substitutions and using (\ref{sympl1}) one sees that the bracket in $\fg$ is as in the statement of Proposition \ref{P1}. Moreover, using that $\w_0$ is a cocycle, one sees that (\ref{jacobi2}), (\ref{jacobi3}), (\ref{jacobi1}) and (\ref{sympl2}) prove that $(D_1,D_2, a_{11},a_{12},a_{21},a_{22})$ is a set of compatible data as defined in Definition~\ref{compd}.
\end{proof}

\section{K\"ahlerian reduction and oxidation}

Since we want to describe K\"ahlerian Lie algebras, we will work with the following definitions:

\begin{definition}{  Let $(\fg,J,\w)$ be a K\"ahler Lie algebra. We say that it admits a {\it K\"ahlerian reduction} if there exists a non-zero isotropic ideal $\fj$ and a vector subspace $V\subset\fg$ complementary to $\fj^\bot$ such that $\fj\oplus V$ is stable under $J$. The reduction will be called {\it central} if $\fj\subset Z(\fg)$. We will say that the reduction is done {\it by a line} or {\it by a plane} whenever the dimension of $\fj$ is, respectively, 1 or 2.
}\end{definition}

 Notice that, in such a case, $\fg=\fj\oplus\fg_0\oplus V$, where $\fg_0=(\fj\oplus V)^\bot$. But $\fg_0=(\fj\oplus V)^\bot$ turns out to be $J$-invariant (because $J$ is Hamiltonian) and non-degenerate for $\w$. Moreover, if  $\w_0=\w_{|\fg_0\times\fg_0}$ and $\overline\w$ is the form induced by $\w$ on $\fj^\bot/\fj$, then the symplectic vector spaces $(\fg_0,\w_0)$ and $(\fj^\bot/\fj,\overline\w)$ are symplectomorphic. We can then define on $\fg_0$ the Lie algebra structure induced by that on $\fj^\bot/\fj$ so that $(\fg_0,\w_0)$ becomes a symplectic Lie algebra. Since $J_0=J_{|\fg_0}$ is obviously Hamiltonian for $\w_0$, we get that $(\fg_0,J_0,\w_0)$ is actually a K\"ahler Lie algebra. 

\begin{definition}{ We will say that $(\fg,J,\w)$ is a {\it K\"ahlerian oxidation} of the K\"ahler Lie algebra  $(\fg_0,J_0,\w_0)$ if $(\fg_0,J_0,\w_0)$ is a K\"ahlerian reduction of $(\fg,J,\w)$. The oxidation will be said to be {\it central, by a line, or by a plane} if  $\fg_0$ is a reduction of the corresponding type.
}\end{definition}

\begin{remark} \label{notared} { It is obvious that every symplectic reduction of a K\"ahler Lie algebra $(\fg,J,\w)$ with respect to a $J$-invariant ideal $\fj$ is actually a K\"ahlerian reduction. In that case, the ideal $\fj$ should obviously be even-dimensional. Nevertheless, we can also have K\"ahlerian reductions with respect to odd-dimensional ideals. For instance, if one chooses an ideal $\fj=\R z$ for some $z\in\fg$ such that $\w(z,Jz)\ne 0$, then the symplectic reduction of $\fg$ with respect to $\fj$ is a K\"ahlerian reduction. 

Furthermore, suppose that $(\fg,J,\w)$ admits a K\"ahlerian reduction with respect to a central isotropic ideal $\fj\subset Z(\fg)$ with $\dim(\fj)\ge 2$ and consider the corresponding decomposition  $\fg=\fj\oplus\fg_0\oplus V$. If $J(\fj)\cap \fj=\{0\}$, as we have that $\fj\oplus V$ is $J$-invariant, there exists a non-zero $z\in\fj$ such that $\omega(z,Jz)\ne 0$ and, as a consequence, $\fg$ admits a central K\"ahlerian reduction  by a line. Further, if we  consider $W=(\fj\oplus V)\cap (\R z\oplus\R Jz)^\bot$, we have that $W$ is a $J$-invariant complement of $\R z\oplus\R Jz$ in $\fj\oplus V$. We can now find a complement $V_1$ to the central isotropic ideal $\fj_1=\fj\cap (\R Jz)^\bot$ in $W$ such that $W=\fj_1\oplus V_1$ and hence $\fg=\R z\oplus\fg_1\oplus \R Jz$, where $\fg_1=(\R z\oplus\R Jz)^\bot=\fj_1\oplus\fg_0\oplus V_1$ admits a K\"ahlerian reduction with respect to $\fj_1$. In the case  $J(\fj)\cap \fj\ne\{0\}$,  one may find $z\in\fj$ such that $Jz\in\fj$ and one has that $\fg$ admits a central K\"ahlerian reduction by the plane $\R z\oplus\R Jz$. The reduced algebra in that case is also K\"ahlerian reducible if $\dim(\fj)> 2$.
This means that central K\"ahlerian reductions can be decomposed as a sequence of successive central K\"ahlerian reductions by lines or by planes. Let us study these two cases in detail.}
\end{remark}

\begin{proposition}\label{kr1}
Let $(\fg,J,\w)$ be a K\"ahlerian Lie algebra and suppose that it admits a K\"ahlerian reduction by a one-dimensional ideal $\fj=\R z\subset Z(\fg)$. Then $(\fg,J,\w)$ is the symplectic oxidation of a K\"ahlerian Lie algebra $(\fg_0,J_0,\omega_0)$ by means of a pair $(D,a)$ where necessarily $a=0$ and $DJ_0=J_0D$.
Conversely, if $(\fg_0,J_0,\omega_0)$ is a Lie algebra with a K\"ahler structure and $D\in \mbox{Der}(\fg_0)$ is such that $[D,J_0]=0$ and  $\cS(D)D$ is Hamiltonian, then the symplectic oxidation of $\fg=\R z\oplus\fg_0\oplus\R v$ by means of $(D,a=0)$ becomes K\"ahlerian with the complex structure defined by $Jx=J_0x$ for all $x\in\fg_0$ and $Jz=\lambda v$, $Jv=-\lambda^{-1}z$, where $\lambda\ne 0$ can be any arbitrary real number. 
\end{proposition}
\begin{proof}
If the reduction of $(\fg,J,\w)$ with respect to  the one-dimensional ideal $\fj=\R z$ is K\"ahlerian, there must be an element $v_0\in\fg$ such that  $\fg=\fj^\bot\oplus\R v_0$ and $\R z\oplus\R v_0$ is $J$-invariant. This implies that $\R z\oplus\R Jz=\R z\oplus\R v_0$ and, hence, $\omega(z,Jz)=\lambda\ne 0$. Let us choose $v=\lambda^{-1}Jz$ so that, for $\fg_0=(\R z\oplus\R v)^\bot$ and according to Proposition~\ref{sde}, we can consider  $\fg=\R z\oplus\fg_0\oplus\R v$ as an oxidation of $\fg_0$ endowed with the K\"ahler structure $(J_0=J_{|\fg_0},\w_0=\w_{\fg_0\times\fg_0})$ by means of a pair $(D,a)\in \mbox{\rm Der}(\fg_0)\times\fg_0$. We have that \[[v,x]=Dx+\omega_0(x,a)z\] for all $x\in\fg_0$ and since $J$ is a complex structure and $Jv\in Z(\fg)$, we have
\begin{eqnarray*}
0&=& [v,x]+J[Jv,x]+J[v,Jx]-[Jv,Jx]
= [v,x]+J[v,Jx]\\ &=&  D x + \omega_0(x,a)z+
J_0 D J_0x +  \omega_0(J_0x,a)Jz,\end{eqnarray*}
for all $x\in\fg_0$. This is only possible if $a=0$ and $DJ_0=J_0D$. 

For the converse, it suffices to see that $J$ is a complex structure and Hamiltonian with respect to the symplectic form $\omega$ in the oxidation. It is obvious that $J^2=-I$. The condition $[D,J_0]=0$ also gives $[D^*,J_0]=0$ and one has that $N_J(z,v)=0$ and
\begin{eqnarray*}
N_J(z,x)&=&\lambda J N_J(v,x)=\lambda J(  D x + \omega_0(x,a)z+
J_0 D J_0x +  \omega_0(J_0x,a)Jz)=0,\\
N_J(x,y)&=&N_{J_0}(x,y)-\w_0(\cS(D)x,y)z-\w_0(\cS(D)J_0x,y)Jz-\w_0(\cS(D)x,J_0y)Jz\\& & +\w_0(\cS(D)J_0x,J_0y)z\\&=& N_{J_0}(x,y)-\w_0(\cS(D)x,y)z+\w_0(J_0\cS(D)x,J_0y)z-\w_0(J_0\cS(D)x,y)Jz\\& &+\w_0(\cS(D)x,J_0y)Jz=0,\end{eqnarray*}
for all $x,y\in\fg_0$. Thus, $J$ is a complex structure. The skew-symmetry of $J$ with respect to $\omega$ is obvious since $\fg_0=(\R z\oplus\R v)^\bot$,  $J_0$ is Hamiltonian with respect to $\w_0$, and one has $\w(v,Jz)=\lambda\w(v,v)=0$, $\w(Jv,z)=\lambda^{-1}\w(z,z)=0$.
\end{proof}

\begin{remark}\label{rema_kahler1} { 
Rescaling $z$ and $v$ in Proposition \ref{kr1} by taking $\tilde z=\sqrt{|\lambda^{-1}|}z$, $\tilde v=\sqrt{|\lambda|}v$, we can consider that $\fg=\R \tilde z\oplus\fg_0\oplus\R \tilde v$ is the corresponding oxidation by means of 
$(\tilde D, a=0)$ where $\tilde D=\sqrt{|\lambda|}D$ and the complex structure on $\fg$ will then be given by $Jx=J_0x$ for $x\in\fg_0$ and $J\tilde z=\varepsilon\tilde v$, where $\varepsilon=\pm 1$ depends on the sign of $\lambda$. 
}\end{remark}

 \begin{proposition} \label{2reduc} Let \( (\fg_0, J_0,  \omega_0) \) be a K\"ahler Lie algebra, and let $(\fg,\omega)$ be the central symplectic oxidation defined by some set of admissible data $(\sH,\sF, (a_{ij})_{i,j\le 2})\in \mbox{\rm Der}(\fg_0)^2\times\fg_0^4$ as in Proposition \ref{P1}. 
If $[J_0\sH-\sF,J_0]=0$ and $J_0(a_{11}-a_{22})=a_{12}+a_{21}$, then the linear map $J:\fg\to\fg$ defined by \[
Jz = z', \quad Jz' = -z, \quad Jv = v', \quad Jv' = -v, \quad Jx = J_0x, \quad \forall x \in \fg_0
\]
is a complex structure on $\fg$ and  \( (\fg, J, \omega) \) is a K\"ahler Lie algebra.
Conversely, if a K\"ahlerian Lie algebra $(\fg, J, \omega)$ admits a central K\"ahlerian reduction $\fg_0$ with respect to a $J$-invariant central plane $\fj$, then  $(\fg, J, \omega)$ is constructed from $(\fg_0,J_0,\w_0)$ that way.
\end{proposition}

\begin{proof}
As we already have that $(\fg,\omega)$ is a symplectic Lie algebra,
we only need to prove that $J$ is a complex structure and that $$\omega(JX,Y)+\omega(X,JY)=0,$$ for all $X,Y \in \fg=\R z\oplus\R z'\oplus\fg_0\oplus\R v'\oplus\R v$. Since one has that $\fg_0=(\R z\oplus\R z'\oplus\R v'\oplus\R v)^\bot$, that both subspaces are stable under $J$, and that $J_0$ is $\w_0$-Hamiltonian, it suffices to see that the restriction of $J$ to  $\R z\oplus\R z'\oplus\R v'\oplus\R v$ is also Hamiltonian. But this is clear because $\R z\oplus\R z'$ and $\R v'\oplus\R v$ are isotropic and for all $\rho,\rho',\xi,\xi'\in\R$ one has:
\begin{eqnarray*}
& & \omega(\rho z+\rho' z',J(\xi' v'+\xi v))=\omega(\rho z+\rho' z',-\xi' v+\xi v')=-\rho \xi' +\rho'\xi,\\
& & \omega(J(\rho z+\rho' z'),\xi' v'+\xi v)=\omega(\rho z'-\rho' z,\xi' v'+\xi v)=\rho \xi' -\rho'\xi.
\end{eqnarray*}
In order to see that $N_J=0$, first notice that one obviously has $N_J(z,\cdot)=N_J(z',\cdot)=0$ because $z,z'=Jz$ are in the center of $\fg$ and that also $N_J(v,v')=0$ because $v'=Jv$. If $x\in\fg_0$, we get
\begin{eqnarray}
N_J(v,x) &=&[v,x]+J[v',x]+J[v,J_0x]-[v',J_0x]\nonumber\\
&=&\sH x+\w_0(x,a_{11})z+\w_0(x,a_{21})z'+J_0\sF x+\w_0(x,a_{12})z'-\w_0(x,a_{22})z\nonumber\\
& & +J_0\sH J_0x+\w_0(J_0x,a_{11})z'-\w_0(J_0x,a_{21})z-\sF J_0x\nonumber\\& &-\w_0(J_0x,a_{12})z-\w_0(J_0x,a_{22})z'\nonumber
\\&=& [J_0 \sH -\sF,J_0](x)+\omega_{0}(J_0x,J_0(a_{11}-a_{22})-a_{12}-a_{21})z\nonumber\\& &+\omega_{0}(x,a_{12}+a_{21}-J_0(a_{11}-a_{22}))z'=0, \label{nijvx}
\end{eqnarray}
because $J_0 \sH - \sF$ commutes with $J_0$ and $J_0(a_{11}-a_{22})=a_{12}+a_{21}$ by hypothesis. This also implies $N_J(v',x)=N_J(Jv,x)=-JN_J(v,x)=0$.
Finally, for $x,y\in\fg_0$ we have
\begin{eqnarray*}
 N_J(x,y) &=&  [x , y] + J[J_0x , y] + J[x, J_0y] - [J_0x , J_0y] \\
&=& N_{J_0}(x,y)-\w_0(\cS(\sH)x,y)z-\w_0(\cS(\sF)x,y)z'\\ & &-\w_0(\cS(\sH)J_0x,y)z'+\w_0(\cS(\sF)J_0x,y)z
 -\w_0(\cS(\sH)x,J_0y)z'\\ & & +\w_0(\cS(\sF)x,J_0y)z+\w_0(\cS(\sH)J_0x,J_0y)z+\w_0(\cS(\sF)J_0x,J_0y)z'\\
        &=& N_{J_0}(x,y) - \omega_{0}(\cS([J_0 \sH -\sF,J_0])  x , y)z  +  \omega_{0}(J_0 \cS([J_0 \sH - \sF,J_0]) x , y)z'\\&=& 0, 
\end{eqnarray*}
because $[J_0 \sH -\sF,J_0]=0$ and $J_0$ is a complex structure.
Thus, the Nijenhuis tensor vanishes identically, which completes the proof.

The converse follows at once from the fact that for any non-zero $z\in\fj$, we have $\fj=\R z\oplus\R Jz$. We can now choose $v\in\fg$ such that $\w(z,v)=1$ and $\w(Jz,v)=0$, so that $\w(Jz,Jv)=1$ and $\w(z,Jv)=0$. Furthermore, we can suppose that $\omega(v,Jv)=0$ by taking $\tilde v=v-\frac{\omega(v,Jv)}{2}Jz$ otherwise. It thus suffices to take $z'=Jz$, $v'=Jv$ and use Proposition~\ref{P2} defining $\alpha=\w([v,v'],v)$, $\beta=\w([v,v'],v')$ and $(\sH,\sF, (a_{ij})_{i,j\le 2})\in \mbox{\rm Der}(\fg_0)^2\times\fg_0^4$ so that the identities 
$$
[v, x] = \sH x+\omega_0(x,a_{11})z + \omega_0(x,a_{21})z',\quad
[v', x] =\sF x+ \omega_0(x,a_{12})z + \omega_0(x,a_{22})z' $$
are verified. We only need to prove that $[J_0 \sH -\sF,J_0]=0$ and $J_0(a_{11}-a_{22})=a_{12}+a_{21}$ but both equalities follow at once from $N_J(v,x)=0$, as in (\ref{nijvx}).
\end{proof}

\section{K\"ahlerian reduction and oxidation of nilpotent Lie algebras}

In this section we will focus on nilpotent Lie algebras admitting K\"ahler structures. The existence of a K\"ahler structure does not impose any  restriction on the nilpotency type of the involved complex structure, and nilpotent Lie algebras with strongly non-nilpotent K\"ahlerian structures are already known \cite{LU}. Moreover, we can give an oxidation method to construct some strongly non-nilpotent K\"ahlerian nilpotent Lie algebras.   

\begin{proposition} \label{SnN} Let $(\fg_0,J_0,\omega_0)$ be a K\"ahlerian nilpotent Lie algebra and let 
\begin{center}$D\in \mbox{\rm Der}(\fg_0)$ be such that $D^2=D^*D=0$, $[D, J_0]=0$, and $\mbox{\rm Ker}(D^*)\cap Z(\fg_0)=\{0\}$,\end{center} and let $(\fg=\R z\oplus\fg_0 \oplus \R v,\omega)$ be the central symplectic oxidation of $\fg_0$ by a line by means of $(D,a=0)$. 
The linear map $J:\fg\to\fg$ defined by $J_{|\fg_0}=J_0$ and $Jz=v$, $Jv=-z$ is a strongly non-nilpotent complex structure and $(\fg,J,\omega)$ is K\"ahler.
\end{proposition}

\begin{proof}
According to Proposition \ref{kr1} and Remark \ref{rema_kahler1}, $(\fg,J,\omega)$ is K\"ahler.  Recall that the non-trivial brackets in $\fg=\R z\oplus\fg_0 \oplus \R v$ are 
$$
[x,y] = [x,y]_0 - \omega_0(\cS(D)x,y)z, \quad
[v,x] =  D x,
$$
for all $x,y\in\fg_0$.

We will see that $Z(\fg)\cap JZ(\fg)=\{0\}$, so that $J$ is strongly non-nilpotent. Suppose that there exist $\mu,\tau\in\R$ and $x_0\in\fg_0$ such that $X=\mu z+x_0+\tau v\in Z(\fg)$ and $$JX=-\tau z+J_0x_0+\mu v\in Z(\fg).$$ Clearly, $0=[X,v]=-Dx_0$ implies that $x_0\in \mbox{\rm Ker}(D)$ and also $J_0x_0\in \mbox{\rm Ker}(D)$ because $DJ_0=J_0D$. Besides, for $y\in\fg_0$ we have
 \begin{eqnarray*}
& & 0=[X,y]=[x_0,y]+\tau Dy=[x_0,y]_0-\w_0(D^*x_0,y)z+\tau Dy,\\
 & & 0=[JX,y]=[J_0x_0,y]+\mu Dy=[J_0x_0,y]_0-\w_0(D^*J_0x_0,y)z+\mu Dy,
\end{eqnarray*} 
from which we deduce $x_0,J_0x_0\in \mbox{\rm Ker}(D^*)$ and $\tau D=-\adgn(x_0)$, $\mu D=-\adgn(J_0x_0)$. 
But then $\adgn(\mu x_0-\tau J_0x_0)=0$, which means that $y_0=\mu x_0-\tau J_0x_0\in \mbox{\rm Ker}(D^*)\cap Z(\fg_0)=\{0\}$, which is only possible if either $x_0=0$ or $\mu=\tau=0$, because $x_0$ and $J_0x_0$ are linearly independent. In this last case, one has $\adgn(x_0)=-\tau D=0$, so that $$x_0\in  \mbox{\rm Ker}(D)\cap\mbox{\rm Ker}(D^*)\cap Z(\fg_0)$$ and we get $x_0=0$ again. Hence, $x_0$ is always the zero vector and  $\mu D=\tau D=0$. But $D$ is not identically zero because  $\mbox{\rm Ker}(D^*)\cap Z(\fg_0)=\{0\}$ and $\fg_0$ has a non-trivial center, and therefore $\mu=\tau=0$, proving $X=0$.
\end{proof}

\begin{example}\label{ejem2} { Let us consider the K\"ahlerian oxidation $\fg_0$ of $\R^2=\mbox{$\R$-span}\{x_1,x_2\}$ with $\omega(x_1,x_2)=1$ and $J_0x_1=x_2$ constructed as in Proposition \ref{2reduc} for $D_1=D_2=0$ and $a_{11}=x_1$, $a_{12}=a_{21}=-x_2$, $a_{22}=3x_1$. It is defined by
the brackets
$$[v,x_1]=-z',\quad [v,x_2]=[v',x_1]=-z,\quad [v',x_2]=-3z',$$
and the symplectic form and the complex structure are defined by
$$\omega(x_1,x_2)=\omega(z,v)=\omega(z',v')=1,\quad Jx_1=x_2,\quad Jv=v',\quad Jz=z'.$$
Let us consider the linear map $D$ on $\fg_0$ defined by $$Dx_1=v,\quad Dx_2=v',\quad Dv=Dv'=Dz=Dz'=0.$$ It is easy to verify that $D$ is a derivation of $\fg_0$.  Obviously, $D^2=0$ and $D^*D=0$ because $D\fg_0\subset\R v\oplus\R v'$ is totally isotropic. Besides, $x\in\mbox{Ker}(D^*)$ if and only if $x\in(\R v\oplus\R v')^\bot$, so that $\mbox{\rm Ker}(D^*)\cap Z(\fg_0)=\{0\}$. According to Proposition \ref{SnN}, we can construct a K\"ahlerian Lie algebra by central symplectic oxidation by means of $(D,a=0)$ with $\alpha=\beta=0$. We then get $\fg=\R z_1\oplus\fg_0\oplus\R v_1$ with the brackets, symplectic form and complex structure given by
\begin{eqnarray}
& & [v_1,x_1]=v,\quad [v_1,x_2]=v',\quad [v,x_1]=-z',\quad [v,x_2]=[v',x_1]=-z,\nonumber\\& &  [v',x_2]=-3z',\quad  [x_1,z]=[x_2,z']=z_1,\label{exa21}\\
& & \omega(x_1,x_2)=\omega(z,v)=\omega(z',v')=\w(z_1,v_1)=1,\label{exa21b}\\& &  Jx_1=x_2,\quad Jv=v',\quad Jz=z',\quad Jz_1=v_1 \label{exa22}.
\end{eqnarray}
 Since  $Z(\fg)$ is one-dimensional, one obviously has that $J$ is a strongly non-nilpotent complex structure. 
}
\end{example}

\begin{example} \label{ejem3} { Let us see that we can oxidize the K\"ahler Lie algebra constructed in Example~\ref{ejem2} to obtain a K\"ahlerian Lie algebra with a weakly non-nilpotent complex structure. Let us rewrite (\ref{exa21}) , (\ref{exa21b}) and (\ref{exa22}) considering $e_0=v_1$, $e_1=v$, $e_2=v'$, $e_3=x_1$, $e_4=x_2$, $e_5=z$, $e_6=z'$ and $e_7=z_1$, so that we will consider the K\"ahlerian Lie algebra $(\fg_0,J_0,\omega_0)$ given by
\begin{eqnarray*}
& & [e_0,e_3]=e_1,\quad [e_0,e_4]=e_2,\quad [e_1,e_3]=-e_6,\quad [e_1,e_4]=[e_2,e_3]=-e_5,\\& & [e_2,e_4]=-3e_6,\quad  [e_3,e_5]=[e_4,e_6]=e_7,\\
& & \omega(e_3,e_4)=\omega(e_5,e_1)=\omega(e_6,e_2)=\w(e_7,e_0)=1,\\ & & Je_1=e_2,\quad Je_3=e_4,\quad Je_5=e_6,\quad Je_7=e_0.
\end{eqnarray*}
Let us take $D_1=D_2=0$ and $a_{11}=a_{22}=-e_7$, $a_{12}=a_{21}=0$. It is clear that $(D_1,D_2, (a_{ij}))$ is a set of compatible data and that the conditions of Proposition \ref{2reduc} are verified. The corresponding K\"ahlerian oxidation for $\alpha=\beta=0$ and that set of data is defined by $\fg=\R \tilde z\oplus\R \tilde z'\oplus\fg_0\oplus\R\tilde v'\oplus\R \tilde v$, where
\begin{eqnarray*}\label{exa2}
& & [\tilde v,e_0]=\tilde z,\quad [\tilde v',e_0]=\tilde z',\quad [e_0,e_3]=e_1,\quad [e_0,e_4]=e_2,\quad [e_1,e_3]=-e_6,\\ &  & [e_1,e_4]=[e_2,e_3]=-e_5,\quad [e_2,e_4]=-3e_6, \quad  [e_3,e_5]=[e_4,e_6]=e_7,\\
& &\omega(\tilde z,\tilde v)=\omega(\tilde z',\tilde v')=\omega(e_3,e_4)=\omega(e_5,e_1)=\omega(e_6,e_2)=\w(e_7,e_0)=1,\\
& & J\tilde z=\tilde z',\quad J\tilde v=\tilde v',\quad Je_1=e_2,\quad Je_3=e_4,\quad Je_5=e_6,\quad Je_7=e_0.
\end{eqnarray*}
One easily verifies that $\fa_1(J)=Z(\fg)\cap JZ(\fg)=\R \tilde z\oplus\R \tilde z'$, and this implies
$$\fa_2(J) =\{ x\in\fg\ \mid \ [x,\fg]\subset \R \tilde z\oplus\R \tilde z', [Jx,\fg]\subset \R \tilde z\oplus\R \tilde z'\}=\mbox{$\R$-span}\{\tilde z,\tilde z',\tilde v,\tilde v'\}.$$
But $\fa_2(J)\cap [\fg,\fg]=\R \tilde z\oplus\R \tilde z'$ and, therefore,
$$\fa_3(J)=\{ x\in\fg\ \mid \ [x,\fg]\subset \fa_2(J), [Jx,\fg]\subset\fa_2(J)\}=\fa_2(J).$$
This clearly shows that $J$ is weakly non-nilpotent.
}
\end{example}

It is clear from Remark \ref{notared} that every K\"ahlerian nilpotent Lie algebra where the complex structure is quasi-nilpotent is K\"ahlerian reducible since $\fa_1(J)=Z(\fg)\cap JZ(\fg)$ is a non-zero $J$-invariant central ideal. The next results give a certain inductive construction of those quasi-nilpotent K\"ahler structures. We start with a straightforward corollary of Proposition \ref{kr1}.

\begin{lemma}  \label{decompcase} Let $(\fg,J,\omega)$ be a K\"ahlerian nilpotent Lie algebra with $\dim(\fg)=n>2$
and suppose that $J$ is quasi-nilpotent. Suppose that there exists $z\in Z(\fg)\cap JZ(\fg)$ such that $\w(z,Jz)\ne 0$, then the K\"ahlerian Lie algebra $(\fg,J,\omega)$ is decomposable and there exists an $(n-2)$-dimensional K\"ahlerian Lie algebra $(\fg_0,J_0,\omega_0)$ such that $\fg=\fg_0\oplus\R^2$, where both ideals are non-degenerate with respect to $\w$ and stable under $J$.
\end{lemma}
\begin{proof}
According to Proposition \ref{kr1}, $\fg$ is a K\"ahlerian oxidation by a line $$\fg=\R z\oplus\fg_0\oplus \R Jz$$ by means of a pair $(D,a=0)$. But $D$ should vanish because $Jz\in Z(\fg)$. The non-trivial brackets then reduce to $[x,y]=[x,y]_0$ for $x,y\in\fg_0$, and the result follows.
\end{proof}

We can now state our main result in this section.  

\begin{theorem}\label{tr:mainth}
Let $(\fg,J,\omega)$ be a K\"ahlerian nilpotent Lie algebra 
endowed with a quasi-nilpotent complex structure and assume 
$\dim (\fg) \ge 2$. Then $(\fg,J,\omega)$ can be obtained by a finite sequence of K\"ahlerian central oxidations by   lines or planes starting from either the Lie algebra $\{0\}$ or a Lie algebra with a strongly non-nilpotent 
K\"ahler structure. Furthermore, all the oxidations by lines in that sequence are actually trivial extensions by a 2-dimensional abelian Lie algebra with a definite K\"ahlerian metric.
If the complex structure $J$ is nilpotent, all the intermediate reductions have a nilpotent complex structure and the 
sequence necessarily finishes at  $\{0\}$.
\end{theorem}
\begin{proof}
Since $J$ is quasi-nilpotent, we can choose a non-zero
$z\in Z(\fg) \cap JZ(\fg).$ If $\omega(z,Jz)$ is non-zero, we can consider the reduction by $\R z$, so that $\fg$ is an oxidation by a line and in fact, according to Lemma \ref{decompcase}, a trivial extension by a 2-dimensional abelian Lie algebra with a definite K\"ahlerian metric. If, on the contrary, $\fj=\R z\oplus\R Jz$ is isotropic, $\fg$ admits a K\"ahlerian reduction to a K\"ahlerian Lie algebra of dimension $(n-4)$ as in Proposition~\ref{2reduc}.

In both cases, let $(\fg_0,J_0,\w_0)$ denote the corresponding reduced algebra. Clearly, if $J_0$ is still quasi-nilpotent, the same argument applies to
$(\fg_0,J_0,\omega_0)$. Repeating this procedure, we obtain
a sequence of K\"ahlerian central reductions by lines or planes. Since the dimension
strictly decreases at each step, the process must terminate.
The iteration stops when one reaches either $\{0\}$  or a Lie algebra carrying a strongly non-nilpotent
K\"ahler structure. Reversing the process, we see that we can obtain $\fg$ by the corresponding successive oxidations.

For the second part of the statement, assume that $J$ is nilpotent and let $(\fg_0,J_0,\w_0)$ be as above. If the oxidation is by a line, according to Lemma \ref{decompcase}, one has that $\fg$ is the direct sum of K\"ahlerian Lie algebras $\fg=\fg_0\oplus\R^2$, so that $\fa_k(J)= \fa_k(J_0)\oplus\R^2$ and, as a consequence, the nilpotency of $J$ clearly implies that of $J_0$.
Now, consider that the oxidation to obtain $\fg$ from $\fg_0$ is done by certain compatible data $(D_1,D_2,(a_{ij}))$. Let us first see by induction on $k$ that $\fg_0\cap\fa_k(J)\subset \fa_k(J_0)$ for every $k\ge 1$. To prove the result for $k=1$, note that if $z_0\in\fg_0\cap Z(\fg)$, then one has for all $x\in\fg_0$ that
$$0=[z_0,x]=[z_0,x]_0-\w_0(\cS(D_1)z_0,x)z-\w_0(\cS(D_2)z_0,x)Jz$$
and this implies $[z_0,x]_0=0$. Therefore,  $\fg_0\cap Z(\fg)\subset Z(\fg_0)$ and,  since $\fg_0$ is stable under $J$, one also has $\fg_0\cap JZ(\fg)\subset J_0Z(\fg_0)$. We then get $\fg_0\cap \fa_1(J)\subset \fa_1(J_0)$.
Let us now suppose that $\fg_0\cap\fa_\ell(J)\subset \fa_\ell(J_0)$ is verified for some $\ell\ge 1$ and let us prove the inclusion for $k=\ell +1$. Since 
$$\fa_{\ell +1}(J)=\{X\in\fg\ \mid \ [X, \fg]\subset \fa_{\ell}(J), [JX, \fg]\subset \fa_{\ell}(J)\},$$
an element $x\in \fg_0$ which belongs to $\fa_{\ell+1}(J)$ in particular satisfies $[x,y],[J_0x,y]\in\fa_{\ell}(J)$ for all $y\in\fg_0$. But one has
\begin{eqnarray*}
& & [x,y]=[x,y]_0-\w_0(\cS(D_1)x,y)z-\w_0(\cS(D_2)x,y)Jz,\\ & &  [J_0x,y]=[J_0x,y]_0-\w_0(\cS(D_1)J_0x,y)z-\w_0(\cS(D_2)J_0x,y)Jz,\end{eqnarray*}
so that  the fact that $z,Jz\in\fa_1(J)\subset \fa_{\ell}(J)$ and $[x,y],[J_0x,y]\in\fa_{\ell}(J)$ clearly imply $[x,y]_0,[J_0x,y]_0\in\fg_0\cap\fa_{\ell}(J)$ and, hence, the induction hypothesis gives $$[x,y]_0,[J_0x,y]_0\in\fa_\ell(J_0).$$ This shows that $x\in\fa_{\ell+1}(J_0)$, which completes the induction. 

Now, since $J$ is nilpotent, there exists $\nu(J)$ such that $\fa_{\nu(J)}(J)=\fg$. We then have
$$\fg_0=\fg_0\cap\fg=\fg_0\cap\fa_{\nu(J)}(J)\subset \fa_{\nu(J)}(J_0).$$
This is only possible if $\fa_{\nu(J)}(J_0)=\fg_0$, so that $J_0$ is nilpotent of step $\nu(J_0)\le \nu(J)$. Therefore, the complex structure of the reduced Lie algebra  must be nilpotent at each reduction. This implies that the reduction cannot start from a strongly non-nilpotent K\"ahler structure, so the process must finish at $\{0\}$.
\end{proof}

Let us now focus on the case of 2-step nilpotent Lie algebras. It is very simple to prove that a complex structure on a 2-step nilpotent Lie algebra is always quasi-nilpotent. Furthermore, it has been seen in \cite[Th. 3]{Gao} that any complex structure $J$ on a non-abelian 2-step nilpotent Lie algebra is nilpotent of step $\nu(J)\in\{ 2,3\}.$ In our next result we characterize the K\"ahlerian oxidations of 2-step nilpotent Lie algebras that give another 2-step nilpotent Lie algebra and give the explicit conditions under which the complex structure on the resulting oxidation has step  $\nu(J)=2$ or $\nu(J)=3$.

\begin{proposition} \label{2-s-nil} Let $(\fg_0,J_0,\w_0)$ be a 2-step nilpotent Lie algebra with a K\"ahlerian structure and consider $(\sH,\sF, (a_{ij})_{i,j\le 2})\in \mbox{\rm Der}(\fg_0)^2\times\fg_0^4$.
The set \[(\sH,\sF, (a_{ij})_{i,j\le 2})\in \mbox{\rm Der}(\fg_0)^2\times\fg_0^4\] is compatible for a K\"ahlerian oxidation and gives a 2-step nilpotent extended Lie algebra if and only if the following conditions hold:
\begin{enumerate}
\item[{\rm (i)}] $D_1,D_2$ are central derivations, $D_iD_j=D_i^*D_j=0$ for $1\le i,j\le 2$, 
 \item[{\rm (ii)}] $a_{ij}\in ([\fg_0,\fg_0]+\mbox{\rm Im}(D_1)+\mbox{\rm Im}(D_2))^\bot$, 
\item[{\rm (iii)}] $a_{12}-a_{21}\in\mbox{\rm Ker}(D_1)\cap\mbox{\rm Ker}(D_2)\cap Z(\fg_0)$,
\item[{\rm (iv)}] $\w_0(a_{12}-a_{21},a_{ij})=0$ for all $i,j\le 2$,
\item[{\rm (v)}] $[J_0D_1-D_2,J_0]=0$, and $J_0(a_{11}-a_{22})=a_{12}+a_{21}$.
\end{enumerate}
In such a case, if $\fg$ is not abelian, then its complex structure $J$ has step $\nu(J)=2$ if and only  if $\nu(J_0)\le 2$, $J_0\mbox{\rm Im}(D_1)\subset {\mathcal W}$, $J_0[\fg_0,\fg_0]\subset {\mathcal W}$ and $J_0(a_{12}-a_{21})\in {\mathcal W}$, where
$${\mathcal W}=\{a\in Z(\fg_0) \ \mid \ D_ia=D_i^*a=0,\ \omega_0(a,a_{ij})=0,\ i,j=1,2\}.$$
\end{proposition}
\begin{proof}
The proof of the first part is quite straightforward. Let us consider the oxidized Lie algebra $\fg=\R z\oplus\R z'\oplus\fg_0\oplus\R v'\oplus\R v$, whose brackets are given as in the statement of Proposition \ref{P1}. Notice that (v) is necessary in order to have a K\"ahlerian oxidation, as in Proposition \ref{2reduc}. From $[v,[v,v']]=[v',[v,v']]=0$, one gets the conditions (iv) and $a_{12}-a_{21}\in\mbox{\rm Ker}(D_1)\cap\mbox{\rm Ker}(D_2)$. When $x,y\in\fg_0$, the identities $[x,[v,y]]=[x,[v',y]]=0$ give that $D_1,D_2$ are central derivations and that $\cS(D_i)D_j=0$ for $i,j\le 2$. But this and condition (i) of Definition \ref{compd} imply that $\adgn(a_{ij})$ is Hamiltonian for all $i,j\le 2$ or, equivalently, since $\w$ is a symplectic form, $a_{ij}\in [\fg_0,\fg_0]^\bot$. The conditions 
$$0=[v,[v,x]]=[v',[v,x]]=[v',[v',x]]=[v,[v',x]]$$ for $x\in\fg_0$ yield $D_iD_j=0$ and $a_{ij}\in (\mbox{\rm Im}(D_1)+\mbox{\rm Im}(D_2))^\bot$. Recalling  that $\cS(D_i)D_j=0$ and combining those equalities with $D_iD_j=0$, one immediately gets $D_i^*D_j=0$ for all $i,j$. Finally, from $[x,[v,v']]=0$, we get the remaining condition $a_{12}-a_{21}\in Z(\fg_0)$. It is not difficult to see that those conditions are also sufficient because they guarantee 
\begin{center}$[v,[x,y]]=[v',[x,y]]=[t,[x,y]]=0$, for all $x,y,t\in\fg_0$\end{center} and conditions (ii) and (iii) of Definition \ref{compd} are also trivially verified. 

Now, let us see under which conditions one gets $\nu(J)=2$. Since, as we showed in the proof of Theorem \ref{tr:mainth}, $\nu(J_0)\le \nu(J)$, we must have $\nu(J_0)\le 2$. Since $\fg$ is not abelian, $\fa_1(J)=Z(\fg)\cap JZ(\fg)\ne \fg$ and, recalling that $\fg$ is 2-step nilpotent, one has $\nu(J)=2$ if and only if  
$$\begin{aligned}
    \fg=\fa_2(J)&=\{X\in\fg\ \mid\ [X,\fg],[JX,\fg]\subset Z(\fg)\cap JZ(\fg)\}\\&=\{X\in\fg\ \mid\ [X,\fg],[JX,\fg]\subset  JZ(\fg)\},
\end{aligned}$$
which is verified if and only if $[\fg,\fg]\subset JZ(\fg)$ or, equivalently, if $J[\fg,\fg]\subset Z(\fg)$. Using the brackets described in the statement of Proposition \ref{P1} and since $z,Jz\in JZ(\fg)$, one easily sees that this occurs if and only if 
\[
    J_0(a_{12}-a_{21})\in Z(\fg), \quad
    J_0D_1x\in Z(\fg), \quad
    J_0D_2x\in Z(\fg) \quad \text{and} \quad
    J[x,y]_0\in Z(\fg),
\]
for all $x,y\in\fg_0$. But one easily sees that
$Z(\fg)\cap \fg_0={\mathcal W}$. Thus, it suffices to show that $J_0D_2x\in {\mathcal W}$ for all $x\in\fg_0$ if $J_0D_1y\in {\mathcal W}$ for all $y\in\fg_0$. To prove it, let us suppose that, for every $y\in\fg_0$, one has $J_0D_1y\in {\mathcal W}$, which means $J_0D_1y\in Z(\fg_0)\cap\mbox{Ker}(D_i)\cap\mbox{Ker}(D_i^*)$  and $\omega_0(J_0D_1y,a_{ij})=0$, for all $i=1,2$. The condition   $[J_0D_1-D_2,J_0]=0$ means
$J_0D_2=D_2J_0-D_1-J_0D_1J_0$ and, using that $D_1,D_2$ are central derivations and that $a_{ij}\in (\mbox{\rm Im}(D_1)+\mbox{\rm Im}(D_2))^\bot$,  we have
\begin{eqnarray*}
& & J_0D_2x=D_2J_0x-D_1x-J_0D_1J_0x\in Z(\fg_0),\\
& & D_iJ_0D_2x=D_iD_2J_0x-D_iD_1x-D_iJ_0D_1(J_0x)=0,\\
& & D_i^*J_0D_2x=D_i^*D_2J_0x-D_i^*D_1x-D_i^*J_0D_1(J_0x)=0,\\
& & \omega_0(J_0D_2x,a_{ij})=\omega_0(D_2J_0x,a_{ij})-\omega_0(D_1x,a_{ij})-\omega_0(J_0D_1(J_0x),a_{ij})=0.\end{eqnarray*}
This completes the proof.
\end{proof}

\begin{remark}{  
In \cite{FPetal} the authors study 2-step nilpotent Lie algebras endowed with a complex structure preserving the center of the algebra. A compact quotient $\Gamma\backslash G$ of the connected and simply-connected Lie group associated to one of those algebras by a uniform discrete subgroup is called a {\it Kodaira manifold}. Inspired by that definition, we will say that a Lie algebra with a complex structure $(\fg_0,J_0)$ is a {\it Kodaira Lie algebra} if  $\fg_0$ is 2-step nilpotent and $JZ(\fg_0)=Z(\fg_0)$.

A Kodaira Lie algebra $(\fg_0,J_0)$  always has a complex structure of step $\nu(J_0)\le 2$ because $\fa_1(J_0)=Z(\fg_0)$ and hence $[\fg_0,\fg_0]\subset \fa_1(J_0)$. Furthermore, for every central derivation $D_1$ one has $J_0D_1x\in J_0Z(\fg_0)=Z(\fg_0)$ and $J_0(a_{12}-a_{21})\in Z(\fg_0)$ whenever $a_{12}-a_{21}\in Z(\fg_0)$. 
However, a 2-step nilpotent K\"ahlerian oxidation of a K\"ahlerian Kodaira Lie algebra is not necessarily a Kodaira Lie algebra. Furthermore, the oxidized Lie algebra can have $\nu(J)=3$. For instance, the 4-dimensional  Lie algebra    $\fg_0={\mathcal K}$ given in Remark \ref{nota3} with the complex structure given by $J_0x_1=x_2$, $J_0x_4=x_3$ is K\"ahlerian. Recall that the bracket and the symplectic form were given by the non-trivial pairings $[x_1,x_2]_0=x_3$, $\w_0(x_1,x_4)=\w_0(x_2,x_3)=1$. It is easy to see that for $D_1=D_2=0$ and $a_{11}=a_{22}=x_1$, $a_{12}=a_{21}=0$, all the conditions of Proposition \ref{2-s-nil} are verified, so that we can construct the corresponding K\"ahlerian oxidation for some arbitrary $\alpha,\beta\in\R$, which will be 2-step nilpotent. According to the explicit expression of the brackets given in Proposition \ref{P1}, one has 
$$
[x_1,x_2]=[x_1,x_2]_0=x_3, \quad [v,Jx_3]=-[v,x_4]=-\w_0(x_4,x_1)z=z.
$$
This means that $x_3\in [\fg,\fg]\subset Z(\fg)$ but $Jx_3\not\in Z(\fg)$, so that $\nu(J)\ne 2$.
}\end{remark}

\end{document}